\documentclass[11pt]{amsart}
\usepackage{latexsym}
\usepackage{amsmath}
\usepackage{amsfonts}
\usepackage{amssymb}
\usepackage{amsthm}
\newtheorem{thm}{Theorem}[section]
\newtheorem{prop}{Proposition}[section] %[thm]
\newtheorem{lem}{Lemma}[section] %[thm]

\newtheorem{Def}{Definition}[section]
\newtheorem{rem}{Remark}[section]

\newenvironment*{Proof}{{\bf Proof.}}

\newcommand{\fa}{\forall}
\newcommand{\dis}{\displaystyle}

\newcommand{\ca}{\mathcal{A}}

\newcommand{\cm}{\mathcal{M}}
\newcommand{\ct}{\mathcal{T}}
\newcommand{\cf}{\mathcal{F}}
\newcommand{\ld}{\ldots}

\newcommand{\el}{\ell}

\newcommand{\ra}{\rightarrow}

\newcommand{\al}{\alpha}
\newcommand{\bi}{\beta}
\newcommand{\ga}{\gamma }
\newcommand{\Ga} {{\varGamma}}

\newcommand{\de}{\delta }

\newcommand{\De} {{\varDelta}}
\newcommand{\e}{\varepsilon }
\newcommand{\f}{\phi}

\newcommand{\thi}{\theta }

\newcommand{\La} {{\varLambda}}

\newcommand{\la}{\lambda }
\newcommand{\mi}{\mu }

\newcommand{\ti}{\tau }
\newcommand{\sm}{\smallsetminus}

\newcommand{\R}{\mathbb{R}}

\newcommand{\Z}{\mathbb{Z}}
\newcommand{\N}{\mathbb{N}}

\newcommand{\qs}{$\quad\square$}

\begin{document}
\title{\bf Sharp weak type inequalities for the dyadic maximal operator}
\author{Eleftherios N. Nikolidakis}
\date{}
\maketitle
\begin{abstract}
We obtain sharp estimates for the localized distribution function of
$\cm\f$, when $\f$ belongs to $L^{p,\infty}$ where $\cm$ is the
dyadic maximal operator. We obtain these estimates given the $L^1$
and $L^q$ norm, $q<p$ and certain weak-$L^p$ conditions.In this way we
refine the known weak (1,1) type inequality for the dyadic maximal
operator.As a consequence we prove that the inequality
\begin{eqnarray}
\|\cm_\ct\f\|_{p,\infty}\le\frac{p}{p-1}||\f||_{p,\infty}
\end{eqnarray}
is sharp allowing every possible value for the $L^1$ and the $L^q$ norm
for a fixed $q$ such that $1<q<p$, where $||\cdot||_{p,\infty}$ is the usual
quasi norm on $L^{p,\infty}$.
\end{abstract}
%
%\vspace{0.3cm}
%
\section{Introduction}  %  1
\noindent

The dyadic maximal operator on $\R^n$ is defined by:
\begin{eqnarray}
\cm_d\f(x)=\sup\bigg\{\frac{1}{|Q|}\int_Q|\f(u)|du:x\in Q, \;
Q\subseteq\R^n \ \ \text{is a dyadic cube}\bigg\} \label{eq1.1}
\end{eqnarray}
for every $\f\in L^1_{loc}(\R^n)$ where $|\cdot|$ is the Lebesgue
measure on $\R^n$ and the dyadic cubes are those formed by the grids
$2^{-N}\Z^n$ for $N=1,2,\ld\;.$

As it is well known it satisfies the following weak type (1,1)
inequality
\begin{eqnarray}
|\{x\in\R^n:\cm_d\f(x)\ge\la\}|\le\frac{1}{\la}\int_{\{\cm_d\f\ge\la\}}|\f(u)|du  \label{eq1.2}
\end{eqnarray}
for every $\f\in L^1(\R^n)$ and every $\la>0$ from which it is easy
to get the following $L^p$ inequality:
\begin{eqnarray}
\|\cm_d\f\|_p\le\frac{p}{p-1}\|\f\|_p.  \label{eq1.3}
\end{eqnarray}
For every $p>1$ and $\f\in L^p(\R^n)$ it is easy to see that the
weak type inequality (\ref{eq1.2}) is best possible and it is proved
in \cite{8} that (\ref{eq1.3}) is also best possible (for general
martingales see \cite{2} and \cite{3}).

In studying the dyadic maximal operator it would be convenient to
work with functions supported in the unit cube $[0,1]^n$ and more
generally defined on a non-atomic probability measure space
$(X,\mi)$ where the dyadic sets are given in a family $\ct$ of
measurable subsets of $X$ that has a tree-like structure similar
to the one in the dyadic case. Then we replace $\cm_d$ by
\begin{eqnarray}
\cm_\ct\f(x)=\sup\bigg\{\frac{1}{\mi(I)}\int_I|\f|d\mi:x\in
I\subseteq X,\; I\in\ct\bigg\}  \label{eq1.4}
\end{eqnarray}
and (\ref{eq1.2}) and (\ref{eq1.3}) remain true and sharp in this
setting.

Actually, in this general setting (\ref{eq1.3}) has been improved
even more by inserting the $L^1$-norm of $\f$ as a variable giving
the so called Bellman functions of the dyadic maximal operator. In
fact in \cite{4} the following function of variables $f,F$ has been
explicitly computed
\begin{eqnarray}
B(f,F)=\sup\bigg\{\int_X(\cm_\ct\f)^pd\mi:\f\ge0, \; \int_X\f
d\mi=f,\; \int_X\f^pd\mi=F\bigg\}  \label{eq1.5}
\end{eqnarray}
where $0<f^p\le F$.

The related Bellman functions for the case $p<1$ have been also
computed in \cite{5}.

It is interesting now to search what happens in case we replace the
$L^p$-norm with the quasi norm $\|\cdot\|_{p.\infty}$ defined in
$L^{p,\infty}$, where
\begin{eqnarray}
\|\f\|_{p,\infty}=\sup\{\la\mi(\{|\f|\ge\la\})^{1/p}:\la>0\}
\label{eq1.6}
\end{eqnarray}
for every $\f$ such that this supremum is finite.

It is known that $L^{p,\infty}$ contains $L^p$ strictly and $\cm_\ct$ can be
defined on $L^{p,\infty}$ with values on $L^{p,\infty}$. As a matter
of fact it is not difficult to see that $\cm_\ct$ satisfies the
following
\begin{eqnarray}
\|\cm_\ct\f\|_{p,\infty}\le\frac{p}{p-1}\|\f\|_{p,\infty}
\label{eq1.7}
\end{eqnarray}
for every $\f\in L^{p,\infty}$.

In \cite{7} it is proved that (\ref{eq1.7}) is best possible.

Actually, a stronger fact is proved there, namely that
\begin{eqnarray}
\sup\bigg\{\|\cm_\ct\f\|_{p,\infty}:\f\ge0,\int_X\f
d\mi=f,\;\|\f\|_{p,\infty}=F\bigg\}=\frac{p}{p-1}F \label{eq1.8}
\end{eqnarray}
for every $(f,F)$ such that $0<f\le\frac{p}{p-1}F$. That is
(\ref{eq1.7}) is sharp allowing every value for the $L^1$-norm of
$\f$.

In the present paper we precisely compute
\begin{equation}
B(f,A,F,\la)=\sup\bigg\{\mi(\{\cm_\ct\f\ge\la\}):\f\ge0, \; \int_X\f
d\mi=f,\; \int_X\f^qd\mi=A, \;  \|\f\|_{p,\infty}=F\bigg\}.
\label{eq1.9}
\end{equation}
for a fixed $q$ such that $1<q<p$, and for all allowable values of
$(f,A,F)$.In order to find (\ref{eq1.9}) it is convenient to work with the following

\begin{equation}
B_{1}(f,A,F,\la)=\sup\bigg\{\mi(\{\cm_\ct\f\ge\la\}):\f\ge0, \; \int_X\f
d\mi=f,\; \int_X\f^qd\mi=A, \;  \|\f\|_{p,\infty}\le F\bigg\}.
\label{eq1.10}
\end{equation}

Our aim is to prove the following
\begin{thm}\label{thm1.1}
If $(f,A,F)\in D_1$, $A>A_2(f,F)$ then $B_1(f,A,F,\la)$ is given by
\[
B_1(f,A,F,\la)=\left\{\begin{array}{ll}
                       1,  & \la\le f \\
                        f/\la, & f<\la\le\frac{p}{p-1}F\la'_3 \\
                        \de', & \frac{p}{p-1}F\la'_3<\la\le\frac{p}{p-1}F\la'_1  \\                  \Big[\frac{\frac{
                        p}{p-1}/F}{\la}\Big]^p, & \la<\frac{p}{p-1}F\la'_1
                      \end{array}\right]
\]
where a) $\la'_1$ is the unique solution of
\[
F_{\la'}\bigg(\frac{1}{(\la')^p}\bigg)+\frac{\Ga}{(\la')^{p-1}}=A
\]
for
\[
\la'=\frac{p-1}{p}\la\frac{1}{F}\ge\bigg(\frac{1}{f_1}\bigg)^{1/(p-1)}, \ \ \text{and} \ \  f_1=\frac{p-1}{p}\cdot\frac{f}{F}.
\]
b) $\la'_3$ is the unique root of the equation
\[
T_{\la'}(f_1/\la')=A \ \ \text{on} \ \ \bigg[f_1,\bigg(\frac{1}{f_1}\bigg)^{1/(p-1)}\bigg].
\]
c) $\de'$ is the unique solution of $T_{\la'}(\de')+F_{\la'}(\de')=A$, for $\de'\in\De$, a suitable interval.

Here $T_\la$ and $F_\la$ are certain functions defined on appropriate intervals, and are defined in the sequel.

Additionally $D_1$ is the domain of the extremal problem defined as follows:
\begin{align*}
(f,A,F)\in D_1\;\Leftrightarrow\,&\text{i)}\;0<f\le\frac{p}{p-1}F \\
&\text{ii)}\;A_2(f,F)\le A\le A_1(f,F)
\end{align*}
where
\[
A_1(f,F)=\bigg(\frac{p-1}{p}\bigg)^{p-q/p-1}\frac{p}{p-q}f^{p-q/p-1}\cdot F^{p(q-1)/(p-1)}
\]
and
\[
A_2(f,F)=\left\{\begin{array}{ll}
                  f^q, & 0<f\le F \\
                  F^q\frac{1}{p-q}\Big\{p-q
                  \Big[p\Big(1-\frac{p-1}{p}\cdot\frac{f}{F}\Big)\Big]^{p-q/p-1}\Big\}, & F<f\le\frac{p}{p-1}F.
                \end{array}\right. \text{\qs}
\]
\end{thm}

As a matter of fact we prove Theorem \ref{thm1.1} in case where $F=\frac{p-1}{p}$. We state it at the beginning of Section 4.

Actually our results generalize the results obtained in \cite{6}.

As an immediate Corollary we obtain
\begin{thm}\label{thm1.2}
The following holds
\begin{align*}
\sup\bigg\{&\|\cm_\ct\phi\|_{p,\infty}:\phi\ge0,\;\int_X\phi d\mi=f,\;\int_X\phi^qd\mi=A,\;\|\phi\|_{p,\infty}=F\bigg\}\\
&=\frac{p}{p-1}F, \ \ \text{for every} \ \ (f,A,F)\in D_1.
\end{align*}
That is (\ref{eq1.8}) is best possible allowing every value of the $L^1$ and $L^q$-norm. \qs
\end{thm}

At last we mention that all the above calculations are independent of the measure space and the associated tree.\ We begin now with:
\section{Preliminaries}  % 2
\noindent

Let $(X,\mi)$ be a non-atomic probability space.

The following holds:
\begin{lem}\label{lem2.1}
Let $\f:X\ra\R^+$ be measurable and $I\subseteq X$ be measurable
with $\mi(I)>0$. Suppose that $\frac{1}{\mi(I)}\int\limits_I\f
d\mi=s$. Then for every $t$ such that $0<t\le\mi(I)$ there exists
a measurable set $E_t\subseteq I$ with $\mi(E_t)=t$ and
$\frac{1}{\mi(E_t)}\int\limits_{E_t}\f d\mi=s$.
\end{lem}
\begin{Proof} Consider the measure space $(I,\mi/I)$ and let
$\psi:I\ra\R^+$ be the restriction of $\f$ on $I$ that is
$\psi=\f/I$. Then if $\psi^\ast:(0,\mi(I)]\ra\R^+$ is the
decreasing rearrangement of $\psi$, we have that
\begin{eqnarray}
\frac{1}{t}\int^t_0\psi^\ast(u)du\ge\frac{1}{\mi(I)}\int^{\mi(I)}_0\psi^\ast
(u)du=s\ge\frac{1}{t}\int^{\mi(I)}_{\mi(I)-t}\psi^\ast(u)du.
\label{eq2.1}
\end{eqnarray}
Since $\psi^\ast$ is decreasing we get the inequalities in
(\ref{eq2.1}), while the equality is obvious since
\[
\int^{\mi(I)}_0\psi^\ast(u)du=\int_I\f d\mi.
\]
From (\ref{eq2.1}) it is easily seen that there exists $r\ge0$
such that $t+r\le\mi(I)$ with
\begin{eqnarray}
\frac{1}{t}\int^{t+r}_r\psi^\ast(u)du=s.  \label{eq2.2}
\end{eqnarray}
It is also easily seen that there exists $E_t$ measurable subset
of $I$ such that
\begin{eqnarray}
\mi(E_t)=t \ \ \text{and} \ \ \int_{E_t}\f
d\mi=\int^{t+r}_r\psi^\ast(u)du  \label{eq2.3}
\end{eqnarray}
since $(X,\mi)$ is non-atomic.

From (\ref{eq2.2}) and (\ref{eq2.3}) we get the conclusion of the
lemma. \qs
\end{Proof}

We now call two measurable subsets $A,B$ of $X$ almost disjoint if
$\mi(A\cap B)=0$.

We give now the following
\begin{Def}\label{def2.1}
A set $\ct$ of measurable subsets of $X$ will be called a tree if
the following conditions are satisfied.
\begin{enumerate}
\item[(i)] $X\in\ct$ and for every $I\in\ct$ we have that
$\mi(I)>0$.
\item[(ii)] For every $I\in\ct$ there corresponds a finite or
countable subset $C(I)\subseteq\ct$ containing at least two
elements such that:
\begin{itemize}
\item[(a)] the elements of $C(I)$ are pairwise almost disjoint
subsets of $I$.
\item[(b)] $I=\cup\, C(I)$.
\end{itemize}
\item[(iii)] $\ct=\bigcup\limits_{m\ge0}\ct_{(m)}$ where
$\ct_0=\{X\}$ and
\[
\ct_{(m+1)}=\bigcup_{I\in\ct_{(m)}}C(I).
\]
\item[(iv)] $\dis\lim_{m\ra+\infty}\sup_{I\in\ct_{(m)}}\mi(I)=0$.
\qs
\end{enumerate}
\end{Def}

From \cite{4} we have the following
\begin{lem}\label{lem2.2}
For every $I\in\ct$ and every $\al$ such that $0<\al<1$ there
exists subfamily $\cf(I)\subseteq \ct$ consisting of pairwise almost
disjoint subsets of $I$ such that
\[
\mi\bigg(\bigcup_{J\in\cf(I)}J\bigg)=\sum_{J\in
\cf(I)}\mi(J)=(1-\al)\mi(I).  \text{\qs}
\]
\end{lem}

Let now $(X,\mi)$ be a non-atomic probability measure space and
$\ct$ a tree as in Definition 2.1. We define the associated maximal
operator to the tree $\ct$ as follows: For every  $\f\in L^1(X,\mi)$
and $x\in X$, then
\[
\cm\f(x)=\cm_\ct\f(x)=\sup\bigg\{\frac{1}{\mi(I)}\int_I|\f|d\mi:\;x\in
I\in\ct\bigg\}.
\]
\section{Domain Of The Extremal Problem} % 3
\noindent

We prove the following:
\begin{thm}\label{thm3.1}
For $f$ and $A$ positive constants the following are equivalent

(i) There exists $\phi:X\ra\R^+$ $\mi$-measurable, non zero such that
\[
\int_X\phi d\mi=f,\ \ \int_X\phi^qd\mi=A,\ \ \|\phi\|_{p,\infty}\le\frac{p-1}{p}.
\]

(ii) $0<f\le1$ and $A_f\le A\le\Ga f^{p-q/p-=1}$ where $\Ga=\Big(\frac{p-1}{p}\Big)^q\cdot\frac{p}{p-q}$ and $A_f$ is defined by:
\[
A_f=\left\{\begin{array}{lll}
            f^q,  & \text{if} & 0<f\le\frac{p-1}{p} \\ [1ex]
             \Big(\frac{p-1}{p}\Big)^q\cdot\frac{1}{p-q}\Big\{p-q[p(1-f)]^{p-q/p-1}\Big\}, & \text{if} & \frac{p-1}{p}<f\le1.
           \end{array}\right.
\]
(We say then that $(f,A)\in D$.).  \qs
\end{thm}
\noindent
{\bf Remark} Actually the domain of the above theorem is the one defined in Theorem \ref{thm1.1}, for $F=\frac{p-1}{p}$.\ A scaling argument on Theorem \ref{thm3.1} then gives the respective domain in Theorem \ref{thm1.1}.

In order to prove the above theorem we will need the following.
\begin{lem}\label{lem3.1}
Let $g_1,g_2:[0,1]\ra\R^+$ be non-increasing functions satisfying the following property:

There exists $c\in(0,1)$ such that
\[
g_1\ge g_2 \ \ \text{on} \ \ (0,c], \ \ \text{and} \ \ g_1\le g_2 \ \ \text{on} \ \ (0,1].
\]
Additionally, $\int\limits^1_0g_1=\int\limits^1_0g_2$.\ Then for every $q>1$ we have that $\int\limits^1_0g^q_1\ge\int^1_0g^q_2$.

If in addition the set $\{g_1>g_2\}\cap(0,c]$ has positive Lesbesgue measure then
\[
\int^1_0g^q_1>\int^1_0g^q_2.
\]
\hfill\text{\qs}
\end{lem}

Lemma \ref{lem3.1} is proved as soon as the following is proved. (A simple approximation argument gives the consequence).
\begin{lem}\label{lem3.2}
Let $\big(a_i\big)^n_{i=1}$, $\big(b_i\big)^n_{i=1}$ be finite $(n\ge3)$, non-increasing sequences of non negative numbers satisfying the following property:

There exists $j\in[1,n)\cap\N$ such that
\[
\begin{array}{lll}
  b_i\ge a_i, & \fa\;i=1,2,\ld,j & \text{and} \\ [1ex]
  b_i\le a_i, & \fa\;i=j+1,\ld,n. &
\end{array}
\]
Additionally $a_1+a_2+\cdots+a_n=b_1+b_2+\cdots+b_n$. Then for every $q>1$ we have that
\begin{eqnarray}
b^q_1+b^q_2+\cdots+b^q_n\ge a^q_1+a^q_2+\cdots a^q_n.  \label{eq3.15}
\end{eqnarray}
If in addition $b_i>a_i$ for some $i\in[1,j]$ then we have strict inequality in (\ref{eq3.15}).  \qs
\end{lem}
\noindent
{\bf Proof of Lemma \ref{lem3.2}}

We prove it by induction on $n\ge3$.\\For $n=3$, we have two sequences $(a_1,a_2,a_3)$ and $(b_1,b_2,b_3)$ such that
\[
a_1\ge a_2\ge a_3\ge0 \ \ \text{and} \ \ b_1\ge b_2\ge b_3\ge0, \ \ \text{with} \ \ a_1+a_2+a_3=b_1+b_2+b_3.
\]
Without loss of generality we also suppose, according to our hypothesis that
\[
b_1\ge a_1, \ \ b_2\le a_2, \ \ b_3\le a_3.
\]
We set $a=a_1+a_2+a_3$ and consider all $(x_1,x_2,x_3)\in\R^3$ such that $x_i\ge0$, $\fa\;i=1,2,3$ and $x_1\ge a_1$, $x_2\le a_2$, $x_3\le a_3$ satisfying the additional condition $x_1+x_2+x_3=a$.

We then set
\[
F(x_1,x_2)=x^q_1+x^q_2+(a-x_1-x_2)^q \ \ \text{with domain}\medskip
\]
\[
S=\left\{\begin{array}{l}
(x_1,x_2):x_1\ge a_1,x_2\le a_2,0\le x_1+x_2\le a,   \\
\hspace*{3cm}x_i\ge0,i=1,2.
\end{array}\right\}
\]
It is not difficult to see (by using the theory of Lagrange multipliers) that: $\dis\min_SF=F(a_1,a_2)$, and that this minimum is attained in only one point, namely $(a_1,a_2)$. So the first step of the induction is easily proved. \medskip

We describe more extensively the inductive step.

We suppose that the Lemma holds for $n\in\N$, $n\ge3$. We prove it for $n+1$.\ So we suppose that
\[
x_1\ge a_1, \ \ x_2\ge a_2,\ld,\ \ x_i\ge a_i,\ \ x_{i+1}\le a_{i+1},\ld,x_n\le a_n,\ \ x_{n+1}\le a_{n+1}
\]
for some $i\in\{1,2,\ld,n\}$, and that
\[
x_1+x_2+\cdots+x_n+x_{n+1}=a_1+a_2+\cdots+a_n+a_{n+1}.
\]

We now set
\[
\begin{array}{lll}
  y_j=x_j, & 1\le j\le n-1 &  \\
  y_n=x_n+x_{n+1} &  & \text{and} \medskip
\end{array}
\]
\[
\begin{array}{ll}
  a'_j=a_j, & 1\le j\le n-1   \\
  a'_n=a_n+a_{n+1}. &
\end{array}
\]
We check now that
\begin{eqnarray}
a'_n\ge y_n \ \ \text{or that} \ \ a_n+a_{n+1}\ge x_n+x_{n+1}.  \label{eq3.16}
\end{eqnarray}

We distinguish two cases

A) If $a_n\ge x_n$ then since $a_{n+1}\ge x_{n+1}$, we have
(\ref{eq3.16}).

B) Let $a_n<x_n$, then by hypothesis we have that $a_{n+1}\ge
x_{n+1}$ and $a_i\le x_i$, $\fa\; i=1,2,\ld,n$. But
$a_1+a_1+\cdots+a_n+a_{n+1}=x_1+x_2+\cdots+x_n+x_{n+1}\ge
a_1+a_2+\cdots+a_{n-1}+x_n+x_{n+1}\Rightarrow$(\ref{eq3.16}). So in any
case (\ref{eq3.16}) is true.

From the induction step we conclude that
\begin{eqnarray}
y^q_1+y^q_2+\cdots+y^q_n\ge(a'_1)^q+(a'_2)^q+\cdots+(a'_n)^q.  \label{eq3.17}
\end{eqnarray}
Then we have:
\begin{align*}
&x^q_1+x^q_2+\cdots+x^q_n+x^q_{n+1}-(a^q_1+a^q_2+\cdots+a^q_n+a_{n+1}^q)\\
=&(y^q_1+y^q_2+\cdots+y^q_{n-1}+y^q_n)-y^q_n+(x^q_n+x^q_{n+1}) \\
&-\big((a'_1)^q+\cdots+(a'_n)^q\big)+(a'_n)^q-a^q_n-a^q_{n+1} \\
=&\big[(y^q_1+y^q_2+\cdots+y^q_n)-\big((a'_1)^q+\cdots+(a'_n)^q\big)\big] \\
&+\big[x^q_n+x^q_{n+1}+(a_n+a_{n+1})^q-\big((x_n+x_{n+1})^q+a^q_n+a^q_{n+1}\big)\big]\\
=&\de_n+\ti_n.
\end{align*}
Here
$\de_n=y^q_1+y^q_2+\cdots+y^q_n-\big((a'_1)^q+\cdots+(a'_n)^q\big)\ge0$
by (\ref{eq3.17}) and so
\begin{align*}
x^q_1+\cdots+x^q_n+x^q_{n+1}-(a^q_1+\cdots+a^q_n+a^q_{n+1})\ge\ti_n=&\big[(a_n+a_{n+1})^q+x^q_n+x^q_{n+1}\big]\\
&-\big[a^q_n+a^q_{n+1}+(x_n+x_{n+1})^q\big].
\end{align*}
It is obvious now according to $(\ref{eq3.16})$ that
$a_n+a_{n+1}\ge x_n\ge x_{n+1}$, and no matter the sequence
$(a_n,a_{n+1},x_n+x_{n+1})$ is arranged with $a_n\ge a_{n+1}$, we
obtain two sequences $(a_n+a_{n+1},x_n,x_{n+1})$,
$(a_n,a_{n+1},x_n+x_{n+1})$ which satisfy the hypothesis of the
lemma in case where $n=3$.\ So we must have that $\ti_n\ge0$.
Additionally if there exists $j\in\{1,2,\ld,\}$ such that $x_j>a_j$
then either $j<n$ so that $\de_n>0$ by the induction step, or $j=n$
and $x_1=a_1,\ld,x_{n-1}=a_{n-1}$. In the latter case we must have
that $x_n>a_n$, $x_{n+1}<a_{n+1}$ $x_n+x_{n+1}=a_n+a_{n+1}$ and we
easily prove that $x^q_n+x^q_{n+1}>a^q_n+a^q_{n+1}$. (This is in
fact the case $n=2$, which is easy to handle).  \qs
\vspace*{0.2cm}

We now prove the following
\begin{thm}\label{thm3.2}
Let $(\ga,\de]$ be a subinterval of $(0,1]$.

Then the following are equivalent
\begin{enumerate}
\item[i)] $\exists\;g:(\ga,\de]\ra\R^+$ non-increasing such that
\[
g(t)\le\psi(t)=\bigg(1-\frac{1}{p}\bigg)t^{-1/p}, \ \ \fa\;t\in(\ga,\de), \ \ \text{and}
\]
\[
\int^\de_\ga g=f_1, \ \ \int^\de_\ga g^q=A_1.
\]
\item[ii)] $0<f_1\le\de^{1-\frac{1}{p}}-\ga^{1-\frac{1}{p}}$ and
$z(f_1)\le A_1\le x(f_1)$, where $z(f_1)$, $x(f_1)$ are defined
below.
\end{enumerate}
\end{thm}
\begin{Proof} i) $\Rightarrow$ ii) Suppose we are given $g$ with the conditions in i). Let
\[
f=\int^\de_\ga\psi=\de^{1-\frac{1}{p}}-\ga^{1-\frac{1}{p}}
\]
\[
A=\int^\de_\ga\psi^q=\Ga\big(\de^{1-\frac{q}{p}}-\ga^{1-\frac{q}{p}}\big), \ \ \text{where}
\]
\[
\Ga=\bigg(\frac{p-1}{p}\bigg)^q\frac{p}{p-q}.
\]
Since $g\le\psi$ on $(\ga,\de]$ we must have
\[
f_1\le f, \ \ A_1\le A.
\]
Consider now the following function defined on $(\ga,\de]$, $g_1(t)=\left\{\begin{array}{cc}
\psi(t),& t\in(\ga,c]\\
0, & t\in(c,\de], \end{array}\right.$, where $c$ is such  that
\begin{eqnarray}
\int^\de_\ga g_1=f_1.   \label{eq3.18}
\end{eqnarray}
Then from (\ref{eq3.18}) we can easily see that $c=\big(\ga^{1-\frac{1}{p}}+f_1\big)^{p/p-1}$. Then
\[
\int^\de_\ga g^q_1=\Ga\Big\{\big(\ga^{1-\frac{1}{p}}+f_1\big)^{p-q/p-1}-\ga^{1-\frac{q}{p}}\Big\}
=:x(f_1).
\]
Given the above $g:(\ga,\de]\ra\R^+$ we easily see that the requirements for Lemma \ref{lem3.1} for the function $g$, $g_1$ are satisfied, so we must have
\[
\int^\de_\ga g^q\le\int^\de_\ga g^q_1=:x(f_1)\Rightarrow A_1\le x(f_1),
\]
that is what we wanted to prove.

We consider now two cases:\vspace*{0.2cm} \\
a) $0<f_1<(\de-\ga)\de^{-\frac{1}{p}}\cdot\frac{p-1}{p}=(\de-\ga)\psi(\de)$.\vspace*{0.1cm} \\
We then must have that
\[
\bigg(\int^\de_\ga g\bigg)^q\le(\de-\ga)^{q-1}\bigg(\int^\de_\ga g^q\bigg),
\]
in view of Holder's inequality, that is $f^q_1\le(\de-\ga)^{q-1}A_1\Rightarrow A_1\ge\frac{f^q_1}{(\de-\ga)^{q-1}}$.\vspace*{0.1cm} \\
b) $(\de-\ga)\frac{p-1}{p}\de^{-1/p}<f_1\le f=\big(\de^{1-\frac{1}{p}}-\ga^{1-\frac{1}{p}}\big)$.\vspace*{0.1cm}

Then we consider the following function $g_2:(\ga,\de]\ra\R^+$ such that
\[
g_2(t)=\left\{\begin{array}{cc}
                (1-\frac{1}{p})c^{-1/p}, & t\in(\ga,c] \\
                \psi(t), & t\in(c,\de],
              \end{array}\right.
\]
with $c$ such that $\int\limits^\de_\ga g_2=f_1$.\ It is easy to see then that $c$ must satisfy:
\begin{eqnarray}
\frac{1}{p}c^{1-\frac{1}{p}}+\ga\frac{p-1}{p}c^{-1/p}=\de^{1-\frac{1}{p}}-f_1.   \label{eq3.19}
\end{eqnarray}
It is, in fact, easy to see that there exists unique $c$ satisfying (\ref{eq3.19}) as soon as $f_1$ is described as in case b) above.\ Then for that $c$ we obtain
\[
\int^\de_\ga g^q_2=(c-\ga)\bigg(\frac{p-1}{p}\bigg)^qc^{-q/p}+\Ga\big(\de^{1-\frac{q}{p}}-
c^{1-\frac{q}{p}}\bigg)=:y(f_1)
\]
where $c:=c(f_1)$ as in (\ref{eq3.19}).

Then using again Lemma \ref{lem3.1} for $g$ and $g_2$ we must have that $A_1=\int\limits^\de_\ga g^q\ge y(f_1)$.\ Defining now
\[
z(f_1)=\left\{\begin{array}{cc}
                f^q_1/(\de-\ga)^{q-1}, & 0<f_1<(\de-\ga)\de^{-1/p}\cdot\frac{p-1}{p} \\
                y(f_1), & (\de-\ga)\frac{p-1}{p}\de^{-1/p}\le f_1\le(\de^{1-\frac{1}{p}}-\ga^{1-\frac{1}{p}})
              \end{array}\right.
\]
we must have that
\[
z(f_1)\le A_1=\int^\de_)\ga g^q\le x(f_1),
\]
that is what we wanted to prove\vspace*{0.1cm} \\
ii) $\Rightarrow$ i) Suppose now that we are given $f_1,A_1$:
\[
\begin{array}{l}
  0<f_1\le\de^{1-\frac{1}{p}}-\ga^{1-\frac{1}{p}} \\ [1ex]
  z(h)\le A_1\le x(f_1).
\end{array}
\]
In the proof of i) $\Rightarrow$ ii) we saw that there exist functions $h_1,h_2:(\ga,\de]\ra\R^+$ non-increasing such that $\int\limits^\de_\ga h_1=\int^\de_\ga h_2=f_1$ and
\[
\int^\de_\ga h^q_1=x(f_1), \quad \int^\de_\ga h^q_2=z(f_1)
\]
with $h_1(t)_1,h_2(t)\le\psi(t)$, for every $t\in(\ga,\de]$.

Consider now the function
\[
T:[0,1]\ra\R^+ \quad \text{such that}
\]
\[
T(\el)=\int^\de_\ga[\el h_1+(1-\el)h_2]^q.
\]
It is obvious that $T$ is continuous on $[0,1]$ and that $T(0)=z(f_1)$ while $T(1)=x(f_1)$.

As a result there exists $\el\in[0,1]$ such that $\int\limits^\de_\ga h^q_\el=A_1$, where $h_\el=\el h_1+(1-\el)h_2$. Additionally $\int^\de_\ga h_\el=f_1$ so $h_\el$ is the wanted function, satisfying the conditions of i).  \qs
\end{Proof}

As a corollary we obtain the following propositions.
\begin{prop}\label{prop3.1}
Let $\al\in(0,1]$ $f_1,A_1$ positive numbers, were $f_1\le
\al^{1-\frac{1}{p}}$. Then the following are equivalent
\begin{enumerate}
\item[(i)] $\exists\;g:(0,\al]\ra\R^+$ Lebesgue measurable such that
$g\le\psi$ on $(0,a]$ and
\[
\int^\al_0g=f_1, \ \ \int^\al_0g^q=A_1.
\]
\item[(ii)] a) If \;$0<f_1\le\dfrac{p-1}{p}\al^{1-\frac{1}{p}}$
\;then\; $\dfrac{f^q_1}{\al^{q-1}}\le A_1\le\Ga
f_1^{p-q/p-1}$\vspace*{0.2cm} \\
 b) If\; $\dfrac{p-1}{p}\al^{1-\frac{1}{p}}\le f_1\le
\al^{1-\frac{1}{p}}$ then
\[
\De_{f_1}(\al)\le A_1\le\Ga f_1^{\frac{p-q}{p-1}}
\]
\end{enumerate}
where
\[
\De_{f_1}(\al)=\bigg(\frac{p-1}{p}\bigg)^p\frac{1}{p-q}\bigg\{p\al^{1-\frac{q}{p}}-q\Big[
p\Big(\al^{1-\frac{1}{p}}-f_1\Big)\Big]\bigg]^{p-q/p-1}\bigg\}.  \text{\qs}
\]
\end{prop}
\begin{prop}\label{prop3.2}
For $a\in(0,1]$ and $f_2,A_2$ such that
$f_2\le1-\al^{1-\frac{1}{p}}$ the following are equivalent
\begin{enumerate}
\item[(i)] $\exists\;g:[a,1]\ra\R^+$ Lebesgue measurable such that
\[
g\le\psi \ \ \text{on} \ \ [\al,1] \ \ \text{and} \ \ \int^1_\al g=f_2,\;\int^1_\al g^q=A_2.
\]
\item[(ii)] a) If \; $f_2\le(1-\al)\dfrac{p-1}{p}$ then
$\dfrac{f^q_2}{(1-\al)^{q-1}}\le A_2\le E_{f_2}(\al)$ where
\[
E_{f_2}(\al)=\Ga\Big[\Big(f_2+\al^{1-\frac{1}{p}}\Big)^{p-q/p-1}-\al^{1-\frac{q}{p}}\Big].
\]
b) If $(1-\al)\dfrac{p-1}{p}\le f_2\le1-\al^{1-\frac{1}{p}}$
then
\[
\Ga_{f_2}(\al)\le A_2\le E_{f_2}(\al)
\]
where
\[
\Ga_{f_2}(\al)=\bigg(\frac{p-1}{p}\bigg)^qc^{-q/p}(c-\al)+\Ga(1-c^{1-\frac{q}{p}})
\]
\end{enumerate}
and $c$ satisfies
\begin{eqnarray}
\frac{1}{p}c^{1-\frac{1}{p}}+\bigg(1-\frac{1}{p}\bigg)\al c^{-1/p}=1-f_2. \text{\qs}  \label{eq4.10}
\end{eqnarray}
\end{prop}
As a consequence, setting $a=1$ in Proposition \ref{prop3.1} we obtain Theorem \ref{thm3.1}.

Using the above arguments one can easily prove that the following is true:
\begin{thm}\label{thm3.2}
For $f$ and $A$ positive constants with $A\neq f^q$ the following are equivalent:
\begin{enumerate}
\item[(i)] There exists $\phi:(X,\mi)\ra\R^+$ measurable, non zero, such that
\[
    \int_X\phi d\mi=f,\quad \int_X\phi^qd\mi=A,  \quad \|\phi\|_{p,\infty}=\frac{p-1}{p}.
\]
\item[(ii)] $0<f\le 1$ and $A_f\le A\le\Ga f^{p-q/p-1}$.

\end{enumerate}
\end{thm}
\section{The Extremal Problem}  % 4
\noindent

We state now Theorem \ref{thm1.1} in case where $F=\frac{p-1}{p}$.
\begin{thm}\label{thm4.1}
Let $(f,A)\in D$, $A>A_f$ and $B_1(f,A,\la)$ given by
\begin{align}
B_1(f,A,\la)=\sup\bigg\{\mi(\{\cm_\ct\phi\ge\la\}):\phi\ge0,\;\int_X\phi d\mi=&f,\; \int_X\phi^qd\mi=f, \nonumber\\
&\|\phi\|_{p,\infty}=\frac{p-1}{p}\bigg\},  \label{eq4.20}
\end{align}
Then the following hold
\[
B_1(f,A,\la)=\left\{\begin{array}{cc}
                      1, & 0<\la<f \\ [1ex]
                      \frac{f}{\la}, & f\le\la<\la_3 \\ [1ex]
                      \de, & \la_3\le\la<\la_1 \\ [1ex]
                      \frac{1}{\la^p}, & \la_1<\la
                    \end{array}\right.
\]
where $\de$ is the unique root of the equation $F_\la(\ga)+T_\la(\ga)=A$, $\la_3$ the unique element of the interval $\Big[f,\Big(\frac{1}{f}\Big)^{1/p-1}\Big]$ such that $T_\la(f/\la)=A$ and $\la_1$ the unique root of the equation $F_\la(1/\la^p)+\frac{\Ga}{\la^{p-q}}=A$ on the interval $\Big[\Big(\frac{1}{f}\Big)^{1/p-1},+\infty\Big)$.

Here $T_\la$, $F_\la$ are certain functions defined on appropriate intervals in the sequel. \qs
\end{thm}

In this section we are going to prove the following
\begin{thm}\label{thm4.2}
If $\al=\al(\la)=B_1(f,A,\la)$ where $(f,A)\in D$ with $A_f<A$, then
\begin{enumerate}
\item[(i)] $\al(\la)=\frac{1}{\la^p}$ for every $\la\ge\la_1$, where $\la_1$ is the unique root of the equation
    \[
    F_\la\bigg(\frac{1}{\la^p}\bigg)+\frac{\Ga}{\la^{p-q}}=A \quad \text{on the interval} \quad \bigg(\bigg(\frac{1}{f}\bigg)^{1/p-1},+\infty\bigg).
\]
\item[(ii)] For every $f<\la<\la_1$, $\al$ equals the supremum of all $\bi$ such that $\bi\le\min\Big\{\frac{f}{\la},\frac{1}{\la^p}\Big\}$ and $f-\la\bi\le1-\bi^{1-\frac{1}{p}}$ for which $T_\la(\bi)\le A-F_\la(\bi)\le S_\la(\bi)$.     \qs
\end{enumerate}
\end{thm}
\begin{Proof} Let $\phi$ be as in (\ref{eq4.20}). We work out the case $\la>\frac{p-1}{p}$.\ A few comments follow for the rest range of $\la$'s.

Fix a $\la>\frac{p-1}{p}$.\ We set $E=\{\cm_\ct\phi\ge\la\}$.

Then $E$ is the almost disjoint union of
elements of $\ct$, $I_j$, $j=1,2,\ld\;.$ Indeed, we just need to
consider those $I\in\ct$ maximal under the condition
$\frac{1}{\mi(I)}\int\limits_I\f d\mi\ge\la$.

For every $j$ we have that
\[
\int_{I_j}\f d\mi\ge\la\mi(I_j).
\]
Summing the above equations up to $j$ we get
\begin{eqnarray}
\int_E\f d\mi\ge\la\mi(E).  \label{eq4.4}
\end{eqnarray}
We again consider the decreasing rearrangement of $\f$, let
$\f^\ast:(0,1]\ra\R^+$. In this point we need a fact which is true
on every non-atomic finite measure space and can be seen in
\cite{1}.

Namely that for every $\de\in[0,1]$
\[
\int^\de_0\f^\ast(t)dt=\sup\left\{\begin{array}{ll}
    & K \; \text{measurable subset} \\ [-2.5ex]
  \int\limits_K\f d\mi: &  \\ [-2.5ex]
   & \text{of}\; X\; \text{such that}\;\mi(K)=\de  \\
\end{array}\right\}
\]
where the supremum is actually attained.

From (4.26) we now get in view of the previous comment for
$a=\mi(E)$ that $\int\limits_0^a\f^\ast(t)dt\ge\la a$, so if we
define by
\[
T(f,A,\la)=\sup\left\{\begin{array}{ll}
                           & \exists\; g:(0,1]\ra\R^+ \ \ \text{decreasing} \\
                          \al\in(0,1]: & \text{such that} \; \int^1_0g=f,\;\int^1_0g^q=A,\;g \le\psi\\
                           & \text{and} \; \int^\al_0g\ge\al\la
                        \end{array}
\right\}
\]
we have that
\begin{eqnarray}
B_1(f,A,\la)\le T(f,A,\la).  \label{eq4.5}
\end{eqnarray}

In fact in relation (\ref{eq4.5}) the converse inequality is also
true. To prove this we need the following.
\end{Proof}
\begin{lem}\label{lem4.1}
If $a\in(0,1]$ and $g:(0,1]\ra\R^+$ are such that
\[
\int^a_0g(t)dt\ge a\la,\;\int^1_0g(t)dt=f,\;\int^1_0g^q(t)dt=A, \ \ g^\ast\le\psi,
\]
then there exists $\f:(X,\mi)\ra\R^+$ measurable with
\[
\int_X\f d\mi=f,\;\int_X\f^qd\mi=A, \;\|\f\|_{p,\infty}\le\frac{p-1}{p}
\]
with the additional property:
\[
\mi(\{\cm\f\ge\la\})\ge a.
\]
\end{lem}
\begin{Proof}
Indeed from Lemma \ref{lem2.2} by setting $I=X$ we guarantee the existence of a sequence $(I_j)_j$ of pairwise almost
disjoint elements of $\ct$ in such a way that
\begin{eqnarray}
\mi(\bigcup_jI_j)=\sum\mi(I_j)=a.  \label{eq4.7}
\end{eqnarray}
We consider the measure space $((0,a],\;|\cdot|)$ where $|\cdot|$ is
Lebesgue measure. Because of $\int\limits^a_0g(t)dt\ge a\la$,
applying Lemma \ref{lem2.1} repeatedly we have as a consequence the
existence of a partition $S=\{A_j,\;j=1,2,\ld\}$ of $(0,a]$, which
consists of Lebesgue measurable subsets of $(0,a]$ such that
\begin{eqnarray}
|A_j|=\mi(I_j) \ \ \text{and} \ \ \int_{A_j}g(t)dt\ge\la|A_j|.  \label{eq4.8}
\end{eqnarray}
For every $j=1,2,\ld$ let now $g_j=(g/A_j)^\ast$ defined on
$(0,|A_j|]$. Since $(X,\mi)$ is non-atomic and $\mi(I_j)=|A_j|$ we
easily see that for every $j$ there exists $\f_j:I_j\ra\R^+$
measurable such that $\f^\ast_j=g_j$. Additionally suppose that
$g'=(g/(a,1])^\ast$ and set $Y=X\sm\cup I_j$. Since $\mi(Y)=1-a$ for
the same reasons we get a $\f':Y\ra\R^+$ such that $(\f')^\ast=g'$.
Then since $I_j$ are pairwise almost disjoint there exists a
measurable function $\f:X\ra\R^+$ such that $\f|_{I_j}=\f_j$ almost
everywhere for every $j$ and $\f|_Y=\f'$. Then it is easy to see
that
\[
\f^\ast=g^\ast\le\psi, \ \  \int_{I_j}\f d\mi=\int_{A_j}gd\mi\ge\la|A_j|=\la\mi(I_j)
\]
where the first equality holds almost everywhere with respect to the Lesbesgue
measure.As a consequence we have that the norm (or quasi norm) condtitions for
$\f$ are satisfied and that
\[
\frac{1}{\mi(I_j)}\int_{I_j}\f d\mi\ge\la \ \ \text{for every} \ \ j=1,2,\ld\;.
\]
So, $\{\cm\f\ge\la\}\supseteq \cup I_j$. As a consequence we get
$\mi(\{\cm\f\ge\la\})\ge a$ and the lemma is\linebreak proved. \qs
\end{Proof}

We continue now the proof of Theorem \ref{thm4.1}. Using Lemma \ref{lem4.1} for $g$ non-increasing we see that
\begin{eqnarray}
B_1(f,A,\la)=T(f,A,\la).   \label{eq4.21}
\end{eqnarray}
In fact we have equality in (\ref{eq4.21}) even if we replace the inequality $\int\limits^a_0g(t)dt\ge\al\la$, given in the definition of $T(f,A,\la)$ by equality, thus getting the function $S(f,A,\la)$ defined by
\[
S(f,A,\la)=sup\left\{\al\in(0,1]:\;\begin{array}{l}
                                     \exists\;g:(0,1]\ra\R^+\;\text{non-increasing} \\  [0.5ex]
                                     \text{such that}\; \int\limits^1_0g=f,\;\int\limits^1_0g^q=A,\;g\le\psi \\ [0.5ex]
                                     \text{and}\;\int\limits^a_0g=\al\la.
                                   \end{array}\right\}
\]
This is true because if $g$ is as in the definition of $T(f,A,\la)$, then of course $\int\limits^a_0g(t)dt\ge\al\la$.\ But then there exists $\bi\ge\al$ such that $\int\limits^\bi_0g(t)dt=\bi\la$ since $\thi(t)$ is a n on-increasing function of $t$, where $\thi$ is defined by $\thi:(0,1]\ra\R^+$ with $\thi(t)=\frac{1}{t}\int\limits^t_0g(u)du$ ($g$ is non-increasing).\ Then we just need to apply Lemma \ref{lem4.1} with $\bi$ in place of $\al$.

From the above we obtain that $B_1(f,A,\la)$ equals the supremum of all $\al\in(0,1]$ for which there exists $g:(0,1]\ra\R^+$ non-increasing and $A_1,A_2\ge0$ such that
\[
\int^a_0g=f,\;\int^a_0g^q=A_1,\;\int^1_\al g=f_2,\;\int^1_\al g^q=A_2
\]
and $g\le\psi$ where $A_1+A_2=A$, $f_1=\la\al$, $f_2=f-\la\al$.

We now take advantage the Proposition \ref{prop3.1} and \ref{prop3.2} of Section 3.

In light of the above mentioned fact and Proposition \ref{prop3.1} for a fixed $\la>f$ and $\la>\frac{p-1}{p}$ we define the following functions
\[
T_\la,S_\la:[0,1/\la^p]\ra\R^+ \ \ \text{by}
\]
\[
T_\la(a)=\left\{\begin{array}{ccc}
                  \la^q\al, & \text{for} & \al
                  \le\Big[\frac{p-1}{p}/\la\Big]^p \\[0.5ex]
                  \De_{f_1}(\al), & \text{for} & \Big[\frac{p-1}{p}/\la\Big]^p<\al\le\frac{1}{\la^p}
                \end{array}\right.
\]
where $f_1=\la\al$ and $S_\la(\al)=\Ga(\la\al)^{(p-q)/(p-1)}$.

In light now of Proposition \ref{prop3.2} we also define:
\[
F_\la,G_\la:\;[0,f/\la]\ra\R^+ \quad\text{for $\al$ such that}\quad f-\la\al\le1-\al^{1-\frac{1}{p}}.
\]

a) If $0<f\le\frac{p-1}{p}$
\[
F_\la(\al)=\frac{(f-\la\al)^q}{(1-\al)^{q-1}} \ \ \text{and} \ \ G_\la(\al)=E_{f_2}(\al)
\]

b) while if $\frac{p-1}{p}<f\le1$
\[
F_\la(\al)=\left\{\begin{array}{ccc}
                    \frac{(f-\la\al)^q}{(1-\al)^{q-1}}, & \text{for} & \frac{f-\frac{p-1}{p}}{\la-\frac{p-1}{p}}\le\al\le\frac{f}{\la} \\ [0.5ex]
                    \Ga_{f_2}(\al), & \text{for} & \al\le\frac{f-\frac{p-1}{p}}{\la-\frac{p-1}{p}}
                  \end{array}\right.
\]
and $G_\la(\al)=E_{f_2}(\al)$, where $f_2=f-\la\al$.

After giving the definitions of $T_\la$, $S_\la$, $F_\la$, $G_\la$ we can say that $B_1(f,A,\la)$ equals the supremum of all $\al\in(0,1]$ such that $\al\le\min\Big\{\frac{f}{\la},\frac{1}{\la^p}\Big\}$ and $f-\la\al\le1-\al^{1-\frac{1}{p}}$ for which there exist $A_1,A_2\ge0$ with
\[
\left.\begin{array}{l}
  T_\la(\al)\le A_1\le S_\la(\al) \\
  F_\la(\al)\le A_2\le G_\la(\al)
\end{array}\right\} \ \ \text{and} \ \ A=A_1+A_2.
\]
We now state the following
\begin{lem}\label{lem4.2}
If $(f,A)\in D$ and $\la>f$ such that
$\al=\al(\la)=B_1(f,A,\la)<\frac{1}{\la^p}$ then there exists
$g:(0,1]\ra\R^+$ such that $g\le\psi$, $\int\limits^\al_0g=\al\la$,
$\int\limits^1_0g=f$, $\int\limits^1_0g^q=A$ and $\int\limits^1_\al
g^q=A_2$ where $A_2=F_\la(\al)$.
\end{lem}
\begin{Proof}
From the definition of $\al(\la)=\al$ and Corollary 4.2 it follows that there exists $g:(0,1]\ra\R^+$ with $g\le\psi$, $\int^1_0g=f$, $\int^\al_0g=\la\al$, $\int^\al_0g^q=A_1$, $\int^1_\al g^q=A_2$, where $A_1+A_2=A$ and $F_\la(\al)\le A_2\le G_\la(\al)$.
In fact, we prove that for every such $g$ we have that $A_2=F_\la(\al)$.

Suppose $g$ is as above with $A_2>F_\la(\al)$. Then $|\{g<\psi\}\cap(0,\al]|\ge\de>0$. This is true since if $g=\psi$ on $(0,\al]$ then
\[
\int^\al_0g=\int^\al_0\psi=\al^{1-\frac{1}{p}}>\la\al, \ \ \text{since} \ \ \al<1/\la^p.
\]
But, then there exists $\bi>\al$ such that $\int^\bi_0g\ge\bi\la$ which gives from the definition of $\al(\la)$ that $\al=\al(\la)\ge\bi$ a contradiction.

Let now $\e>0$ be small enough (we will arrange it's choice at the end), with $A_2-\e>F_\la(\al)$.

Since $F_\la(\al)<A_2-\e<G_\la(\al)$, by Proposition 4.2 there exists $g_1:[\al,1]\ra\R^+$ such that $g_1\le\psi$, $\int^1_\al g_1=f-\la\al$, $\int^1_\al g^q_2=A_2-\e$. We set $g_1=g$ on $(0,\al]$.

Since now $|\{g_1<\psi\}\cap(0,\al]|\ge\de>0$ we can construct a function $g_2:(0,1]\ra\R^+$ with the same integral as $g_1$ and larger $L^q$-norm.This is done by increasing it's values on $(0,\al]$ and decreasing in appropriate way it's values on $[\al,1]$, so that $g_2\le\psi$. If $\e>0$ small enough we can arrange everything (since $\de$ is fixed positive) in such a way that $\int^1_0g^q_2=A$. Additionally we have that $\int^1_0g_2=f$ and $\int^\al_0g_2>\la\al$, because we have increased the values of $g_1$ to $g_2$ in the interval $(0,\al]$. This gives as before the existence of a $\bi>\al$ such that $\int^\bi_0g_2\ge\la\bi$, so $\al(\la)\ge\bi>\al$, a contradiction. \qs
\end{Proof}

In the sequel we state and prove the following
\begin{lem}\label{lem4.3}
For $(f,A)\in D$ such that $ A>\ca_f$ and
$\al(\la)=B_1(f,A,\la)$, there exists
$\la_1\ge\Big(\frac{1}{f}\Big)^{1/p-1}$ such that
$\al(\la)=\frac{1}{\la^p}$, for every $\la\ge\la_1$.
\end{lem}
\begin{Proof} If $\la\ge\Big(\frac{1}{f}\Big)^{1/p-1}$ then
$\frac{1}{\la^p}\le\frac{f}{\la}$.

We consider the equation
\[
F_\la(1/\la^p)=A-\frac{\Ga}{\la^{p-q}}.
\]
We easily see that
\[
\lim_{\la\ra\infty}F_\la(1/\la^p)=\ca_f< A=\lim_{\la\ra+\infty}\bigg(A-\frac{\Ga}{\la^{p-q}}\bigg).
\]
For $\la=\la_0=\Big(\frac{1}{f}\Big)^{1/p-1}$ we have that
\[
F_{\la_0}(1/\la^p_0)=F_{\la_0}\bigg(\frac{f}{\la_0}\bigg)=0\ge A-\Ga f^{p-q/p-1}=A-\frac{\Ga}{\la^{p-q}_0}.
\]
So, there exists $\la\ge\la_0$ such that
\[
F_\la(1/\la^p)=A-\frac{\Ga}{\la^{p-q}}.
\]
Let
\[
\la_1=\inf\bigg\{\la\ge\la_0:F_\la(1/\la^p)=A-\frac{\Ga}{\la^{p-q}}\bigg\}
\]
which is obviously a minimum. Then
\begin{eqnarray}
F_{\la_1}\bigg(\frac{1}{\la^p_1}\bigg)=A-\frac{\Ga}{\la_1^{p-q}}.  \label{eq4.11}
\end{eqnarray}
Consider the interval $\big[\frac{1}{\la^p_1},1\big]$. Applying Proposition
\ref{prop3.2} for $\al=\frac{1}{\la^p_1}$,
$f_2=f-\frac{1}{\la^{p-1}_1}$ we obtain that there exists
$g_2:\Big[\frac{1}{\la_1^p},1\Big]\ra\R^+$ such that
\[
g_2\le\psi\bigg/_{\Big[\frac{1}{\la_1^p},1\Big]}, \; \int^1_{1/\la_1^p}g_2=f-\frac{1}{\la_1^{p-1}},\;
\int^1_{1/\la_1^p}g_2^q=F_{\la_1}\bigg(\frac{1}{\la^p_1}\bigg).
\]
But then if $g:(0,1]\ra\R^+$ with
$g\Big/_{\Big(0,\frac{1}{\la_1^p}\Big]}=\psi$ and
$g\Big/_{\Big[\frac{1}{\la_1^p},1\Big]}=g_2$ we have because of
the above that
\[
\int^1_0g_1=f,\;\int^1_0g^q=A,\;g\le\psi,\;\int^{1/\la^p_1}_0g=\frac{1}{\la_1^{p-1}}=\frac{1}
{\la_1^p}\cdot\la_1
\]
and according to Lemma \ref{lem4.1} we have that
$\al(\la_1)=B_1(f,A,\la_1)\ge\frac{1}{\la^p_1}$. But of course
$\al(\la_1)\le\frac{1}{\la^p_1}$, so that
$\al(\la_1)=\frac{1}{\la_1^p}$.Then we easily see that
$\al(\la)=\frac{1}{\la^p}$ for every $\la\ge\la_1$. This is true
because $g:(0,1]\ra\R^+$ as mentioned before satisfies:
\begin{eqnarray}
g\le\psi,\;\int^{1/\la^p}_0g=\frac{1}{\la^{p-1}}=\frac{1}{\la^p}\cdot\la  \label{eq4.12}
\end{eqnarray}
for every such $\la$ and
$\int\limits^1_0g=f$, $\int\limits^1_0g^q=A$,for every such $\la$.  \qs
\end{Proof}

Let now $\la_2=\min\Big\{\la:\al(\la)=\frac{1}{\la^p}\Big\}$ and $\la$ such that $\al(\la)=\frac{1}{\la^p}$.\ Then $\frac{1}{\la^p}\le\frac{f}{\la}\Rightarrow\la\ge\Big(\frac{1}{f}\Big)^{1/p-1}=\la_0$, so that
\[
\la_2=\min\bigg\{\la\ge\la_0:\al(\la)=\frac{1}{\la^p}\bigg\}.
\]
Let now $\la_1$ as defined in Lemma \ref{lem4.3}.

Obviously $\la_1\ge\la_2$.

Since $\la_2$ is the minimum positive
$\la$ such that $\al(\la)=1/\la^p$, we have from Lemma \ref{lem4.3}
and by continuity reasons that there exists $g:(0,1]\ra\R^+$ such
that
\[
\int^{1/\la^p_2}_0g=\frac{1}{\la^{p-1}_2}=\frac{1}{\la^p_2}\cdot\la_2,\;
\int^1_0g=f,\;\int\limits^1_0g^q=A,\;g\le\psi
\]
and such that
\[
A_2=\int^1_{1/\la^p_2}g^q=F_{\la_2}(\al)
\]
where $\al=\al(\la_2)=\frac{1}{\la^p_2}$. But then
\[
T_{\la_2}(\al(\la_2))=S_{\la_2}(\al(\la_2))=\frac{\Ga}{\la_2^{p-q}}.
\]
Since
\[
T_{\la_2}(\al)\le A_1=\int^{1/\la^p_2}_0g^q\le S_{\la_2}(\al)
\]
we have that
\[
F_{\la_2}\bigg(\frac{1}{\la^p_2}\bigg)+\frac{\Ga}{\la_2^{p-q}}=A
\]
that is $\la_1=\la_2$.

Let now $\la=\la_1=\la_2$.

Then
\[
F_\la(1/\la^p)+\frac{\Ga}{\la^{p-q}}=A,
\]
and $\la\ge\Big(\frac{1}{f}\Big)^{1/p-1}$. For $\mi>\la$ and
$\bi=\frac{1}{\mi^p}$ we have that
\[
\frac{1}{\mi^p}\le\frac{f}{\mi} \ \ \text{and} \ \ f-\mi\bi\le1-\bi^{1-\frac{1}{p}}.
\]
Then $F_\mi(\bi)=F_\mi(1/\mi^p)$ describes the minimum $L^q$-norm
value of functions $g$ defined on
\[
\bigg[\frac{1}{\mi^p},1\bigg]=[\bi,1] \ \ \text{for which} \ \ \int^1_\bi g=f-\mi\bi \ \ \text{and} \ \ g\le\psi.
\]
So
\[
F_\mi\bigg(\frac{1}{\mi^p}\bigg)+\frac{\Ga}{\mi^{p-q}}=\int^1_0g^q_\mi
\]
where $g_\mi$ is defined such that
\[
g_\mi:=\psi,[0,\bi],\ \ \int^1_\bi g^q_\mi=F_\mi\bigg(\frac{1}{\mi^p}\bigg)
\ \ \text{and} \ \ g_\mi\le\psi.
\]
Since $q>1$ we have that $\mi\mapsto\int^1_0 g^q_\mi$ is decreasing on $\big(\big(\frac{1}{f}\big)^{1/(p-1)},+\infty\big)$.

We provide a proof for this right now.

We distinguish two cases\vspace*{0.1cm}

(i) $0<f\le\frac{p-1}{p}\Rightarrow 0<f-\frac{1}{\mi^{p-1}}\le\frac{p-1}{p}$, $ \fa\;\mi\in\Big(\Big(\frac{1}{f}\Big)^{1/(p-1)},+\infty\Big)$.\vspace*{0.1cm}

Then $g_\mi$ has the form
\[
g_\mi(t)=\left\{\begin{array}{cc}
                  \psi(t), & t\in\Big(0,\frac{1}{\mi^p}\Big] \\ [2ex]
                  c_\mi, & t\in\Big[\frac{1}{\mi^p},1\Big],
                \end{array}\right.\ \ \text{where}
\]
\[
c_\mi=\frac{f-\frac{1}{\mi^{p-1}}}{1-\frac{1}{\mi^p}}, \ \ \text{while}
\ \ \int^1_{1/\mi^p}g^q_\mi=F_\mi(1/\mi^p).
\]
It is now easy to see that $v>\mi>\Big(\frac{1}{f}\Big)^{1/(p-1)}\Rightarrow c_v>c_\mi$.\ Additionally $\int\limits^1_0g_v=\int\limits^1_0g_\mi$.

According now to the way that $g_v$ and $g_\mi$ are defined we obtain, by using Lemma \ref{lem3.1} that $\int\limits^1_0g^q_\mi>\int\limits^1_0g^q_v$.\ That is $\mi\mapsto\int\limits^1_0g^q_\mi$ is strictly decreasing on the range $\Big(\Big(\frac{1}{f}\Big)^{1/(p-1)},+\infty\Big)$ for the case $0<f<\frac{p-1}{p}$. \vspace*{0.1cm}

(ii) We consider now the second case $\frac{p-1}{p}<f\le1$.\vspace*{0.1cm}

Let $\mi_0$ be such that $f-\frac{1}{\mi^{p-1}_0}=\frac{p-1}{p}$, so that
\[
\mi_0=\bigg[\frac{1}{\Big(f-\frac{p-1}{p}\Big)}\bigg]^{1/(p-1)}>\bigg(\frac{1}{f}\bigg)^{1/(p-1)}.
\]
Then (A): for every $\Big(\frac{1}{f}\Big)^{1/(p-1)}<\mi\le\mi_0\Rightarrow0<f-\frac{1}{\mi^{p-1}}
\le\frac{p-1}{p}$:  while (B): for $\mi\ge\mi_0\Rightarrow f-\frac{1}{\mi^{p-1}}\ge\mi_0$:\;(B).

Case (A) is worked out as in case (i), and so we see that $\mi\mapsto\int\limits^1_0g^q_\mi$ is strictly decreasing on the range $\Big(\Big(\frac{1}{f}\Big)^{1/(p-1)},\mi_0\Big]$.

Let $\mi$ be as in (B), that is $\mi>\mi_0$. Then we define
\[
g_\mi(t)=\left\{\begin{array}{cc}
                  \psi(t), & t\in\Big(0,\frac{1}{\mi^p}\Big] \\ [2ex]
                  c'_\mi, & t\in\Big(\frac{1}{\mi^p},c_\mi\Big] \\ [2ex]
                  \psi(t), & t\in(c_\mi,1]
                \end{array}\right.
\]
where $c_\mi$ is defined (see Proposition \ref{prop3.2}) by:
\begin{eqnarray}
\frac{1}{p}c^{1-\frac{1}{p}}_\mi+\bigg(1-\frac{1}{p}\bigg)c^{-1/p}_\mi=1-f_2=1-
\bigg(f-\frac{1}{\mi^{p-1}}\bigg)  \label{eq4.22}
\end{eqnarray}
and $c'_\mi=\psi(c_\mi)=\Big(1-\frac{1}{p}\Big)c^{-1/p}_\mi$. From (\ref{eq4.22}) we have that
\begin{eqnarray}
\frac{1}{p}c^{1-\frac{1}{p}}_\mi+\bigg(1-\frac{1}{p}\bigg)c^{-1/p}_\mi-\frac{1}
{\mi^{p-1}}=1-f.  \label{eq4.23}
\end{eqnarray}
It is not now difficult to see, by differentiating (\ref{eq4.23}) with variable $\mi$ that the function $C:(\mi_0,+\infty)\ra\R$ defined by $C(\mi)=c_\mi$ is strictly decreasing, so according again to Lemma \ref{lem3.1} we have that $v>\mi\Rightarrow\int\limits^1_0g^q_\mi>\int\limits^1_0g^q_v$, that is $\mi\mapsto\int\limits^1_0g^q_\mi$ is strictly decreasing on $(\mi_0,+\infty)$. But it is easy to see that $\mi\mapsto\int^1_0g^q_\mi$ is continuous on $\Big(\Big(\frac{1}{f}\Big)^{1/p-1},+\infty\Big)$ so that our assertion is proved.

From the above facts we obtain that for every $\mi>\la=\la_1=\la_2$ we have that
\[
F_\mi\bigg(\frac{1}{\mi^p}\bigg)+\frac{\Ga}{\mi^{p-q}}<A.
\]
The conclusion is Theorem \ref{thm4.2}.
\section{Proof Of Theorem \ref{thm4.1}}
\noindent

We continue now with the proof of {\bf Theorem \ref{thm4.1}}\vspace*{0.2cm}\\
\noindent
{\bf Proof of Theorem 4.1.} We will use Theorem \ref{thm4.1} (i) and (ii).

Let us consider the case $f<\la<\la_1$, so that $\al=\al(\la)=B_1(f,A,\la)<\frac{1}{\la^p}$.\ Of course, we must also have that $\al\le\frac{f}{\la}$.\  We search now for those
$\bi\in\Big[0,\frac{1}{\la^p}\Big]$ such that
$f-\la\bi\le1-\bi^{1-1/p}$.

Consider $K$ defined on $\Big[0,\frac{1}{\la^p}\Big]$ by
$K(\bi)=f-1+\bi^{1-1/p}-\la\bi$. Since
$K'(\bi)=\frac{p-1}{p}\bi^{-1/p}-\la$, $K$ increasing on
$[0,\bi_0]$, decreasing on $\Big[\bi_0,\frac{1}{\la^p}\Big]$ with
maximum value at the point $\bi_0$ where
$\bi_0=\Big[\frac{(p-1)/p}{\la}\Big]^p$. Then
\[
K(\bi_0)=f-1+\bigg[\frac{(p-1)/p}{\la}\bigg]^{p-1}\cdot\frac{1}{p}
\]
which may be positive as well as negative. We first work in case
that $K(\bi_0)>0$ and $\frac{p-1}{p}<f\le1$. From the above we have
that there exist $\bi_1,\bi_2\le\frac{1}{\la^p}$ with $\bi_1<\bi_2$
so that $f-\la\bi_i=1-\bi_i^{1-\frac{1}{p}}$ for $i=1,2$ and for
$\bi\le\frac{1}{\la^p}$ we have that
$f-\la\bi\le1-\bi^{1-\frac{1}{p}}$ if and only if
$\bi\in[0,\bi_1]\cup\Big[\bi_2,\frac{1}{\la^p}\Big]$. With the above
hypothesis we prove the following
\begin{lem}\label{lem4.4}
For $(f,A)\in D$ such that $A>\ca_f$, $f<\la<\la_1$ we have that
\[
\al=B_1(f,A,\la)\in\bigg[\bi_2,\min\bigg\{\frac{f}{\la},\frac{1}{\la^p}\bigg\}\bigg].
\]
\end{lem}
\begin{Proof}
Obviously, for $\ga=\frac{f}{\la}$,
$f-\ga\la\le1-\ga^{1-\frac{1}{p}}$ and
$\la\bi_2\le\bi_2^{1-\frac{1}{p}}$ so by means of Proposition
\ref{prop3.1} there exists $\f:[0,\bi_2]\ra\R^+$ such that
\[
\int^{\bi_2}_0\f=\la\bi_2,\;\f\le\psi,\;\int^{\bi_2}_0\f^q=T_\la(\bi_2).
\]
Now since $\bi_2>\bi_0=\Big[\frac{p-1/p}{\la}\Big]^p$,
$T_\la(\bi_2)=\De_f(\bi_2)$. We extend now $\f$ on $[0,1]$ by
defining $\f=\psi$ on $[\bi_2,1]$. Then since
$f-\la\bi_2=1-\bi_2^{1-\frac{1}{p}}$ we have that
$\int\limits^1_0\f=f$.

By definition now of $\f$ and $\De_f(\bi_2)$ the form of $\f$ must
be such that $\int\limits^1_0\f^q=\ca_f$. This is true since $\De_f(\bi_2)$ describes the $L^q$-norm of a function of the form
\[
h(t):=\left\{\begin{array}{ll}
               \frac{p-1}{p}c^{-1/p}, & t\in(0,c] \\ [0.5ex]
               \psi, & t\in[c,1]
             \end{array}\right.
\]
Since then $\phi=\psi$ on $[\bi_2,1]$ we should have that
\[
\phi(t):=\left\{\begin{array}{lcc}
             \frac{p-1}{p}c^{-1/p}, & \text{on} & (0,c] \\ [0.5ex]
             \psi, & \text{on} & [c,1]
           \end{array}\right.
\]
almost everywhere and $\int^1_0\phi=f$ that is $\int^1_0\phi^q=\ca_f$, by the definition of $\ca_f$ (see Lemma 3.2).
But we remind that
$\int\limits^{\bi_2}_0\f=\la\bi_2$. Since now $\int^1_0\phi^q=\ca_f<A$, and $(f,A)\in D$ we can construct a function $g:(0,1]\ra\R^+$ such that $g\le\psi$, $\int^1_0g=f$, $\int^1_0g^q=A$ (by increasing suitably of $\phi$ on an interval of the form $(0,c_1]$, $c_1\le c$ and decreasing it on $[c_1,1]$) and $\int^\ga_0g\ge\ga\la$, for a $\ga>\bi_2$. But then $B_1(f,A,\la)=\al\ge\ga>\bi_2$. Obviously then $\al\in\big[\bi_2,\min\big\{\frac{f}{\la},\frac{1}{\la^p}\big\}\big]$. \qs
\end{Proof}

Consider now the following function defined on
\[
R_\la:\De=\bigg[\bi_2,\min\bigg\{\frac{f}{\la},\frac{1}{\la^p}\bigg\}\bigg]\ra\R^+
\]
with $R_\la(\bi)=F_\la(\bi)+S_\la(\bi)$, $\bi\in\De$.

We prove now that $R_\la$ is increasing on $\De$. We are in the case where $\frac{p-1}{p}<f\le1$, and $R_\la$ defined on $\De=\Big[\bi_2,\min\Big\{\frac{f}{\la},\frac{1}{\la^p}\Big\}\Big]$.

We have that
\[
R_\la(\ga)=F_\la(\ga)+S_\la(\ga), \ \ \text{for} \ \ \ga\in\De.
\]
Here $F_\la(\ga)=\int\limits^1_\ga g^q_\ga$, where $g_\ga$ is defined as
\[
g_\ga(t)=\left\{\begin{array}{ll}
                  \Big(1-\frac{1}{p}\Big)c^{-1/p}_2, & t\in(0,c_2] \\ [1ex]
                  \psi(t), & t\in(c_2,\ga] \\ [1ex]
                  \Big(1-\frac{1}{p}\Big)c_1^{-1/p}, & t\in(\ga,c_1] \\[1ex]
                  \psi(t), & t\in(c_1,1],
                \end{array}\right.
\]
where $c_1,c_2$ are such that
\[
\int^1_\ga g_\ga=f-\la\ga \ \ \text{and} \ \ \int^\ga_0g_\ga=\la\ga.
\]
Then
\[
S_\la(\ga)=\int^\ga_0g^q_\ga, \ \ \text{and} \ \ R_\la(\ga)=\int^1_0g^q_\ga.
\]
Easily we see that $c_2(\ga)c_2=[p\ga(\ga^{-1/p}-\la)]^{p/(p-1)}$ and that $\ga_1>\ga$, $\ga_1,\ga_2\in\De\Rightarrow c_2(\ga_1)<c_2(\ga)$.

Additionally $c_1=c_1(\ga)$ satisfies
\[
\La(\ga)=\frac{1}{p}[c_1(\ga)]^{1-\frac{1}{p}}+\frac{p-1}{p}\ga[c_1(\ga)]^{-1/p}
-\la_\ga=1-f.
\]
Differentiating with respect to $\ga$ we obtain
\begin{eqnarray}
c'_1(\ga)\frac{1}{p}\cdot\frac{p-1}{p}[c_1(\ga)]^{-1/p}[1-(c_1(\ga))^{-1}\ga]=\la+
\frac{p-1}{p}c_1(\ga)^{-1/p}.  \label{eq5.24}
\end{eqnarray}
But $c_1(\ga)>\ga$, for every $\ga\in\De$.

So from (\ref{eq5.24}) we have that $(c_1(\ga))$ is an increasing function of $\ga$, $\ga\in\De$. Considering again Lemma \ref{lem3.1} we see that $\ga<\ga_1\Rightarrow\int\limits^1_0g_{\ga_1}^q>\int\limits^1_0 g^q_\ga$ or that $R_\la(\ga_1)>R_\la(\ga)$.

The other cases are treated in similar ways.

Then we have that $R_\la(\ga)$, $\ga\in\De$ is an increasing function.

Let now $\al=B_1(f,A,\la)$ then as we mentioned before, $\al\in\De$,
and of course $T_\la(\al)\le A-F_\la(\al)\le S_\la(\al)$ because of
Lemma \ref{lem4.3}. If
$\al<\min\Big\{\frac{f}{\la},\frac{1}{\la^p}\Big\}$ and
$F_\la(\al)+T_\la(\al)<A$ then since
$f-\la\al<1-\al^{1-\frac{1}{p}}$ there exists $\ga$ such that
$\al<\ga<\min\Big\{\frac{f}{\la},\frac{1}{\la^p}\Big\}$ such that
$F_\la(\ga)+T_\la(\ga)<A$ and of course
\[
R_\la(\ga)=F_\la(\ga)+S_\la(\ga)\ge A \ \ (R_\la\;\text{is increasing)}.
\]
That is
\[
T_\la(\ga)\le A-F_\la(\ga)\le S_\la(\ga).
\]
Then Theorem \ref{thm4.2} gives $B_1(f,A,\la)\ge\ga>\al$, a
contradiction that is if
$\al=B(f,A,\la)<\min\Big\{\frac{f}{\la},\frac{1}{\la^p}\Big\}$ we
must have that $F_\la(\al)+T_\la(\al)=A$.

Consider now $\la_0=\Big(\frac{1}{f}\Big)^{1/p-1}$ and the function
$h:E=\Big[f,\Big(\frac{1}{f}\Big)^{1/(p-1)}\Big]\ra\R^+$ defined by
$h(\la)=T_\la(f/\la)$. Notice that for
$f\le\la\le\Big(\frac{1}{f}\Big)^{1/(p-1)}$ we have that
$\frac{f}{\la}\le\frac{1}{\la^p}$, so this definition makes sense. Notice also that $\la_1\ge\la_0$.

Then
\[
h(f)=T_f(1)=\ca_f< A \ \ \text{and} \ \ h(\la_0)=\Ga\frac{1}{\la^{p-q}_0}=
\Ga f^{p-q/p-1}\ge A.
\]
Consider now all $\la\in E$ such that
\[
h(\la)=T_\la(f/\la)\le A.
\]
Then
\[
T_\la(f/\la)\le A=A-F_\la(f/\la)\le\Ga f^{p-q/p-1}=S_\la(f/\la)
\]
so in view of Corollary 4.2, $B_1(f,A,\la)=\frac{f}{\la}$.

For $\la\in E\cap\{\la:h(\la)>A\}$ we have by the above comments that
\[
B_1(f,A,\la)=\sup\{\al\in\De:\;F_\la(\al)+T_\la(\al)=A\}.
\]
Additionally the same equality holds for $\la$ such that $\big(\frac{1}{f}\big)^{1/p-1}<\la\le\la_1$.

So we found $B_1(f,A,\la)$ as $f/\la$ or a maximum root of an equation in case where $\la\in[f,\la_1]$, $K(\bi_0)>0$, $\frac{p-1}{p}<f\le1$.

The case $K(\bi_0)=0$ is worked out in the same way where we replace
$\bi_2$ by $\bi_0$, while the case $K(\bi_0)<0$ is worked out for
\[
\De=\bigg[0,\min\bigg\{\frac{f}{\la},\frac{1}{\la^p}\bigg\}\bigg].
\]
Analogous results are obtained when $0<f\le\frac{p-1}{p}$ where
\[
\De=\bigg[0,\min\bigg\{\frac{f}{\la},\frac{1}{\la^p}\bigg\}\bigg],
\]
since then
\[
f-\la\bi\le\frac{p-1}{p}(1-\bi)\le1-\bi^{1-\frac{1}{p}} \ \ \text{for every}
\ \ \bi\le\frac{1}{\la^p}.
\]

We prove now the following fact.

The function $h(\la)=T_\la(f/\la)$ defined on $E=\Big[f,\Big(\frac{1}{f}\Big)^{1/(p-1)}\Big]$ as above is strictly increasing.

We consider again two cases.\vspace*{0.1cm}

i) $0<f<\frac{p-1}{p}$.\ Then for $f<\la<\Big(\frac{1}{f}\Big)^{1/(p-1)}\cdot\Big(\frac{p-1}{p}\Big)^{p/(p-1)}$ we have that $\frac{f}{\la}\le\Big[\frac{(p-1)/p}{\la}\Big]^p$ so that $h(\la)=\la^{q-1}f$, which is certainly strictly increasing on this range.

As for the range $\Big(\frac{1}{f}\Big)^{1/(p-1)}\Big(\frac{p-1}{p}\Big)^{p/(p-1)}<\la<\Big(\frac{1}{f}\Big)^{1/p-1}$ we have that $\Big[\frac{(p-1)/p}{\la}\Big]^p\le\frac{f}{\la}\le\frac{1}{\la^p}$, so that
\[
h(\la)=\De_{\la_\al}(\al)=\int^\al_0g^q_\la, \ \ \text{where}
\]
\[
g_\la(t)=\left\{\begin{array}{cc}
                  \frac{p-1}{p}c^{-1/p}_2, & t\in(0,c_2] \\ [1ex]
                  \psi(t), & t\in(c_2,\al]
                \end{array}\right.
\]
where $c_2<\al$ is such that $\int\limits^\al_0g_\la=\la\al=\la\frac{f}{\la}=f$.

We prove that $h(\la)$ is increasing on this range also.\ Since it is continuous on all $E$ we obtain the desired result. Easily
\[
c_2(\la)=c_2=\bigg\{p\bigg[\bigg(\frac{f}{\la}\bigg)^{1-\frac{1}{p}}-\la\al\bigg]\bigg\}^{p/p-1}.
\]
So $c_2(\la)$ is decreasing on the above range and thus by using Lemma \ref{lem3.1} we obtain what we wanted to prove.\vspace*{0.1cm}

ii) The second case $\frac{p-1}{p}<f\le1$ is easier to handle because $h(\la)=T_\la(f/\la)$ is defined in a certain way in all $E$, as we have done in the case i) (second part).

We need to say a few comments for the case $f<\la<\frac{p-1}{p}$.

In this case we have the same results with $F_\la(\al)$, $T_\la(\al)$ given by:
\[
F_\la(\al)=\frac{(f-\la\al)^q}{(1-\al)^{q-1}} \ \ \text{and} \ \ T_\la(\al)=\la^q\al
\]
and then
\[
B_1(f,A,\la)=\sup\bigg\{\al\in\De:\frac{f-\la\al)^q}{(1-\al)^{q-1}}+\la^q\al=A\bigg\}
\]
for $\De=\Big[\bi'_2,\min\Big\{\frac{f}{\la},\frac{1}{\la^p}\Big\}\Big]$ for some suitable constant $\bi'_2$.

This finishes the proof of Theorem \ref{thm4.1}.
\begin{rem}\label{rem4.1}

i) The case where $A=\ca_f$ can be worked out separately because
there exists essentially a unique function $g:[0,1]\ra\R^+$ such that
$\int\limits^1_0g=f$, $\int\limits^1_0g^q=A$, $g\le\psi$.

ii) We have that $B(f,A,\la)=B_1(f,A,\la)$, for $A\neq f^q$ as
mentioned in the beginning of this section. This is true of course
for $\la\ge\la_1$, that is for $\la$ such that
$\al(\la)=B_1(f,A,\la)=\frac{1}{\la^p}$.

Now for $\la<\la_1$ let $\al=B_1(f,A,\la)$. Then there exists
$g:[0,1]\ra\R^+$ such that $\int\limits^1_0g=f$,
$\int\limits^1_0g^q=A$, $\int\limits^\al_0g=\al\la$, $g\le\psi$.
Then it is easy to see that for every $\e>0$ small enough we can
change $g$ to $g_\e$ in a way that
\[
\int^{\al-\e}_0g_\e\ge(\al-\e)\la, \ \ \int^1_0g_\e=f, \ \ \int^1_0g^q_\e=A+\de_\e,
\ \ \|g_\e\|_{p,\infty}=\frac{p-1}{p}
\]
and $\de_\e\ra0$ as $\e\ra0^+$. This using continuity arguments
gives $B(f,A,\la)=\al$.

iii) Notice the continuity of the function as calculated on Theorem
\ref{thm4.2}, at the point $\la=\la_1$. As a matter of fact $\de$ is
such that $F_{\la_1}(\de)+T_{\la_1}(\de)=A$. But $\la_1$ is such
that
\[
F_{\la_1}\bigg(\frac{1}{\la_1^p}\bigg)+\frac{\Ga}{\la_1^{p-q}}=A, \ \ \text{and}
\ \ \frac{\Ga}{\la^{p-q}_1}=S_{\la_1}\bigg(\frac{1}{\la_1^p}\bigg)=T_{\la_1}
\bigg(\frac{1}{\la^p_1}\bigg).
\]
So that
\[
F_{\la_1}(\de)+T_{\la_1}(\de)=F_{\la_1}\bigg(\frac{1}{\la^p_1}\bigg)+T_{\la_1}
\bigg(\frac{1}{\la_1^p}\bigg)=A,
\]
which in view of Remark iii) above, gives $\de=\frac{1}{\la^p_1}$.
\end{rem}

Theorem \ref{thm1.2} is now an immediate consequence.

\end{document}